\documentclass{llncs}

\usepackage{llncsdoc}
\usepackage[table]{xcolor}
\usepackage{booktabs}
\usepackage{amsmath}
\usepackage{amsfonts}
\usepackage{amssymb}
\usepackage{latexsym}
\usepackage{exscale,cmmib57}
\usepackage{color}
\usepackage{multirow}
\usepackage{algorithm}  
\usepackage{algorithmic}
\usepackage{epic}
\usepackage{array}
\usepackage{epic,eepic}
\usepackage{fancybox}
\usepackage{psfrag}
\usepackage{caption}
\usepackage{float}
\usepackage{multirow}
\usepackage{graphicx}
\usepackage{epstopdf}
\usepackage{caption}
\usepackage{subfig}
\DeclareMathOperator\erf{erf}
\def\beqns{\begin{eqnarray*}}
\def\eeqns{\end{eqnarray*}}

\newcommand*{\vcenteredhbox}[1]{\begingroup
\setbox0=\hbox{#1}\parbox{\wd0}{\box0}\endgroup}
\makeatletter
\def\blfootnote{\gdef\@thefnmark{}\@footnotetext}
\makeatother

\begin{document}

\title{Quasi-interpolation on a sparse grid \\ with Gaussian}

\author{
Fuat Usta\inst{1} \and  Jeremy Levesley \inst{2} }

\institute{Department of Mathematics, Duzce University , \\
Konuralp Campus, 81620, Duzce, Turkey, \\
\email{fuatusta@duzce.edu.tr},\\ 
\and
Department of Mathematics, University of Leicester, \\
University Road, Leicester, LE1 7RH, United Kingdom, \\
\email{j.levesley@le.ac.uk}
}
\blfootnote{\textup{2010} \textit{Mathematics Subject Classification}:
65F10, 65N55, 65N22, 65D32 }

\maketitle

\begin{abstract}
Motivated by the recent multilevel sparse kernel-based interpolation (MuSIK) algorithm proposed in [Georgoulis, Levesley and Subhan, SIAM J. Sci. Comput., 35(2), pp. A815-A831, 2013], we introduce the new quasi-multilevel sparse interpolation with kernels (Q-MuSIK) via the combination technique. The Q-MuSIK scheme achieves better convergence and run time in comparison with classical quasi-interpolation; namely, the Q-MuSIK algorithm is generally superior to the MuSIK methods in terms of run time in particular in high-dimensional interpolation problems, since there is no need to solve large algebraic systems.

We subsequently propose a fast, low complexity, high-dimensional quadrature formula based on Q-MuSIK interpolation of the integrand. We present the results of numerical experimentation for both interpolation and quadrature in high dimension.
\end{abstract}

\section{Introduction}
Over the last half century, numerical methods have obtained much attention, not only among mathematicians but also in the scientific and engineering communities.

High dimensions usually cause some problems for mathematical modelling on gridded data. The main problem is known as the curse of dimensionality, a term due to Bellmann \cite{Bell}. There is a exponential relationship between the computational cost of approximation with a given error bound $ \epsilon $ and the dimension $ d $ of the space $ \mathbb{R}^{d} $ for any given problem. For this reason, classical approximation techniques are limited to low dimensions. For example, the complexity of solving an approximation problem on a gridded data over a bounded domain $ \Omega \in \mathbb{R}^{d} $ is $ \mathcal{O}(N^{d}) $, where $ N $ is the dimension of the input data.

Radial Basis Function (RBF) interpolation is one of the tools which is effective in approximating and interpolating high-dimensional functions when the data is scattered in its domain. The use of RBFs has become increasingly popular as an approximation method since it obtains delicate and accurate consequences without using a mesh. This is true not only in approximation or interpolation of data sets \cite{powell} but also in solving partial differential equations; see \cite{kansa2}, for example of applications of the RBF method.

Another powerful tool for multidimensional problems is quasi-interpolation, which is comprehensively used in scientific computations, mechanics and engineering, and from which a number of successful results have been gained \cite{beatsonpow,Wu2,Wu}. A quasi-interpolation technique based on radial basis functions is discussed in the sequel. One of the advantages of quasi-interpolation is relevant to features of the generating functions themselves, such as smoothness, simplicity, good shape properties and their exponential decay behaviour at infinity. In addition to these, quasi-interpolation can yield a solution directly since there is no need to solve any large-scale system of equations. Therefore this method can approximate the function in a reduced computational time, even in high dimension, in comparison with other meshless techniques such as RBF interpolation. Quasi-interpolation has been successfully applied to scattered data approximation and interpolation, numerical solutions of partial differential equations and quadrature. 

In the literature of the quasi-interpolation, some methods having convergence rates have been also discussed \cite{mazya1,mazya2}. Although quasi-interpolation operators have principally used functions defined all over the real line, some applications require functions defined on a specific compact interval to allow for efficient approximation procedures, such as boundary integral equations and treatment of partial differential equations. Müller and Varnhorn \cite{muller} applied quasi-interpolation operators to such functions. In contrast to all-space functions, a truncation error has to be controlled for these kind of functions; in addition to these pointwise error estimates and $L_1$ error estimates have also to be given explicitly \cite{muller}. 

Another study of quasi-interpolation has been made by Chen and Cao \cite{Cao1}. In this study, the convergence analysis in the supremum norm has been presented by modifying the quasi-interpolation operator. Then, in \cite{Cao2}, further investigation of quasi-interpolation has been made on the compact interval.

In order to cope with high dimensional problem, there are a number of remedies proposed in academia. The first of these is hyperbolic cross spaces, which were mentioned by Babenko \cite{Babenko} and Smolyak \cite{Smolyak} in the structure of numerical integration in early 1960's. In these papers, the hyperbolic cross product was proposed in order to construct a quadrature rule for multidimensional functions via additional smoothness assumptions. Then, sparse grid methods were introduced by Zenger \cite{Zenger} in 1991 as a solution for PDEs. These techniques have also been used for approximation and interpolation. It can be seen that these methods are similar to hyperbolic cross product hypothesis. The $ sparse $ $ grid $ $ method $ arises from a sparse tensor product construction (hyperbolic cross product) and hierarchical basis. Additionally, sparse grid techniques have also been used as solutions of PDEs for finite difference and finite volume methods. 

In 2012, Georgoulis, Levesley and Subhan \cite{fazli2} introduced a new kernel-based interpolation technique which circumvents both computational complication and conditioning problems. They obtained more reliable and faster results in higher dimension RBF interpolation problems by means of this technique, which can be seen as an improved version of hyperbolic cross functions. The basic principle of this method is the use of anisotropic radial basis function interpolation. All in all, sparse interpolation with kernels, called SIK, allows for an enormous reduction in the amount of computing resources required to solve interpolation problems at a given level of accuracy.

In order to take the advantages of the SIK method a step further, we have introduced the quasi-sparse interpolation with kernels, called Q-SIK. The main motivation for the proposed method stems from the idea of the SIK method. The principal point of the proposed algorithm is to use anisotropic Gaussian interpolation for all directions associated with the sparse grid combination technique. The combination technique \cite{Schneider} is a sparse grid representation, where partial solutions are computed on a certain sequence of coarser grids; one then gets the solution by linearly combining all partial solutions. This is the main idea behind the proposed algorithm. 

In order to obtain both accelerated convergence and more accurate consequences for the interpolation, multilevel techniques have been suggested by a number of researchers. The first to consider the multilevel approximation were Floater and Iske \cite{Iske}. They combined a thinning algorithm and compactly supported radial basis function interpolation. Furthermore, multilevel interpolation techniques enable us to combine the benefits of stationary and non-stationary standard RBFs interpolation, such that this leads to an accelerated convergence. Other researchers have since contributed the multilevel interpolation literature \cite{Hales}.

The proposed scheme, Q-SIK, can be extended to a multilevel variant, quasi-multilevel sparse interpolation with kernels (Q-MuSIK). The Q-MuSIK technique is derived from the Q-SIK method. The main considerations behind the approach to multilevel techniques can be divided into two steps: the first is to interpolate the data sets at the coarsest level, and the second is to update the interpolation of the residuals on gradually finer data sets by means of properly scaled basis functions. 

Principally, construction of a sparse grid in Q-SIK techniques is similar to multilevel methods. When we interpolate given data sets using the Q-MuSIK method, we need to use sparse grids which are nested, and from lower to higher levels, i.e., $ \mathcal{S}_{n,d} \subset \mathcal{S}_{n+1,d}$ where $d$ is dimension and $ n \in \mathbb{N} $. In addition to this, properly scaled anisotropic quasi-interpolation should be used for each level. All in all, it can be said that the Q-MuSIK method does not adversely affect the complexity features of SIK because of the geometric progression in the problem dimension in the nested-ness of sparse grids.

The new method of multilevel sparse kernel-based interpolation (MuSIK) was introduced by Georgoulis, Levesley and Subhan \cite{fazli2} in 2013. This method uses anisotropic radial basis function interpolation on partial grids and then linearly combines them with the combination technique (SIK method). In other words, SIK can be obtained by solving many sub-interpolation problems with properly selected sub-grids and then linearly combining them. Then by combining SIK and multilevel algorithm, the new approach MuSIK has been obtained. In other words the MuSIK uses a residual interpolation at each levels. This technique provides us with a considerable advantage in the interpolation of huge amounts of data in multiple dimensions.

In this context, the MuSIK has been extended in \cite{DGLF} both
theoretically and numerically. More detailed this study showed that SIK and accordingly MuSIK scheme are interpolatory for the special case of scaled Gaussian kernels. In addition to this a numerical integration algorithm is also proposed in \cite{DGLF}, based on interpolating the (high-dimensional) integrand. A series of numerical examples have been presented, highlighting the practical applicability of the proposed algorithm for both interpolation and quadrature for up to 10-dimensional problems in this study.

For a detailed discussion on MuSIK technique, its theoretical explanation and numerical results we refer the reader to \cite{DGLF}, and references therein.
 
\section{Sparse quasi-interpolation with kernels}
One of the most powerful aspects of quasi-interpolation is that it can be generalized easily to the multidimensional case because of its simple structure. In more detail, we have seen some basis functions which have a number of nice properties for quasi-interpolation on a uniform grid in $\mathbb{R}^n$ such as simplicity, smoothness and exponential decay behaviour at infinity. Thus, this kind of basis function makes multidimensional approximation possible with basic analytical representations.

Another advantage of quasi-interpolation is that there is no need to solve large algebraic systems. For example, for some  approximation techniques, such as radial basis functions interpolation, one needs to solve large number of matrix systems to find the best approximation and require significant storage for univariate problems solutions. However, due to the nature of quasi-interpolation, we just need to know function value at specific points. Thus, this makes quasi-interpolation useful in practice not only for low dimensional approximations but also in multidimensional cases. Also integration using quasi interpolation methods provides us positive weights because of its nature. 

Now, in order to provide more benefits from nice properties of quasi-interpolation, we will propose to apply this new technique to quasi-interpolation, which is called the quasi-sparse kernel-based interpolation (Q-SIK) method. Thus, we similarly expect to obtain results with lower computational cost, especially in the interpolation of huge amount of data in multiple dimensions.

Let $ \Omega:=[0,1]^d $, $ d\geq 2 $, and let $ u: \Omega \rightarrow \mathbb{R}$. In order to clarify sparse grid construction, we need to introduce some multi-index notation. Throughout this and next sections we will use $\textbf{l}:=(l_1,\ldots,l_d)\in \mathbb{N}^d$ and $\textbf{i}:=(i_1,\ldots,i_d)\in \mathbb{N}^d$ as a spatial position. The key point here is that the multi-index $\textbf{l}$ consists of non zero components. In other words, non-zero elements are an undesired condition for solving interpolation problems.

For the above multi-indices, the family of directionally uniform grids $\Upsilon_\textbf{l}$ on $\Omega$ can be described with mesh size $h_\textbf{l};=2^{-\textbf{l}}=(2^{-l_1},\ldots,2^{-l_d})$. That is, the family of grids $\Upsilon_\textbf{l}$ consists of the points
\begin{equation}
\textbf{x}_{\textbf{l},\textbf{i}}:=(x_{l_1,i_1},\ldots,x_{l_d,i_d}),
\end{equation}
where $x_{l_p,i_p}:=i_p \cdot h_{l_p}=i_p \cdot 2^{-l_p}$ and $i_p \in \lbrace 0,1,\ldots,2^{l_p} \rbrace$. These grids may have different mesh sizes for each coordinate direction. In addition to this, one can calculate the number of nodes $ \mathcal{P}_{\textbf{l}}$ in $\Upsilon_\textbf{l}$ by using the formula
\begin{equation}
\mathcal{P}_{\textbf{l}}:=\prod_{p=1}^{d}(2^{l_p}+1).
\end{equation}
For instance, should one choose constant $\textbf{l}$, which is equal to $n$ for all $p=1,\ldots,d$, the Cartesian full grid converts to a uniform full grid denoted by $\Upsilon^{n,d}$ with $\mathcal{P}=(2^n+1)^d$, where $n$ is the grid level. We also consider the following subset of $\Upsilon^{n,d}$,
\begin{equation}
\Upsilon_{n,d}:=\bigcup_{|\textbf{l}|_1=n+d-1}\Upsilon_\textbf{l}
\end{equation}
with $ |\textbf{l}|_1:=l_1+l_2+\cdots+l_d $, which will be referred to as the sparse grid of level $  n$ in $  d$ dimensions. We refer to Figure 1 for a visual representation of (3) for $n = 4$ and $d = 2$. Notice that there is some redundancy in this definition as some grid points are included in more than one sub-grid.
\begin{figure}
\begin{center}
\begin{tabular}{m{1.8cm}m{0.2cm}m{1.6cm}m{0.2cm}m{1.6cm}m{0.2cm}m{1.6cm}m{0.2cm}m{1.6cm}}
   \includegraphics[width=1.9cm]{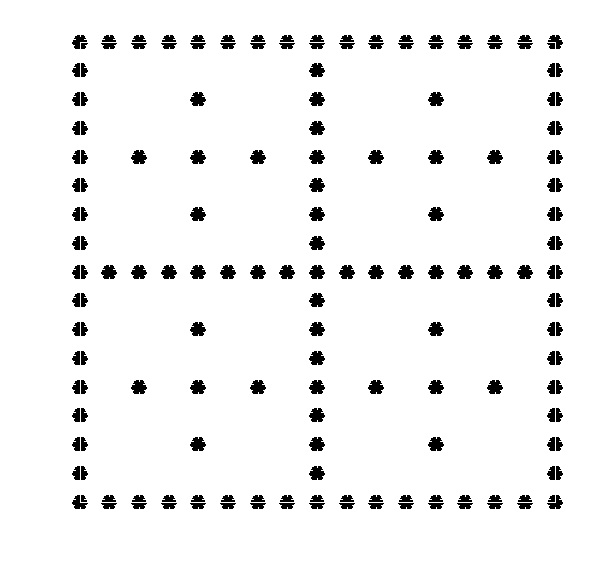} &=
      &\includegraphics[width=1.7cm]{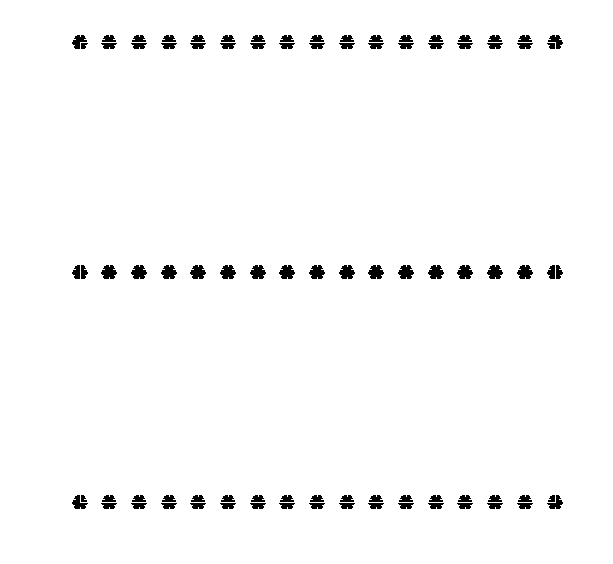} 
     & $\ \cup$  &\includegraphics[width=1.7cm]{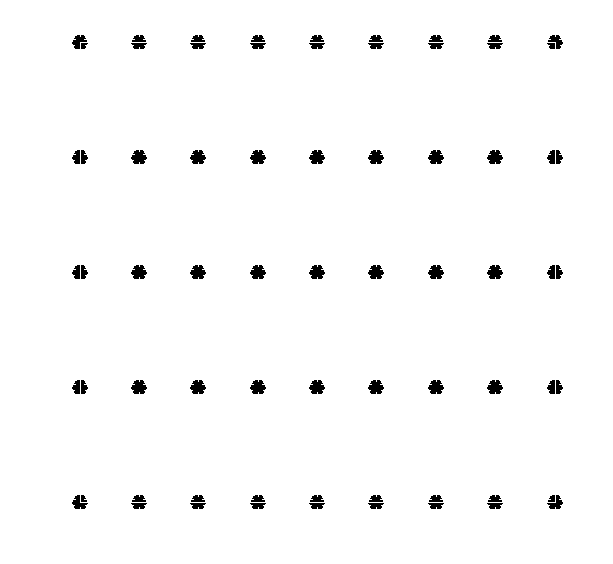} 
     & $\ \cup$ & \includegraphics[width=1.7cm]{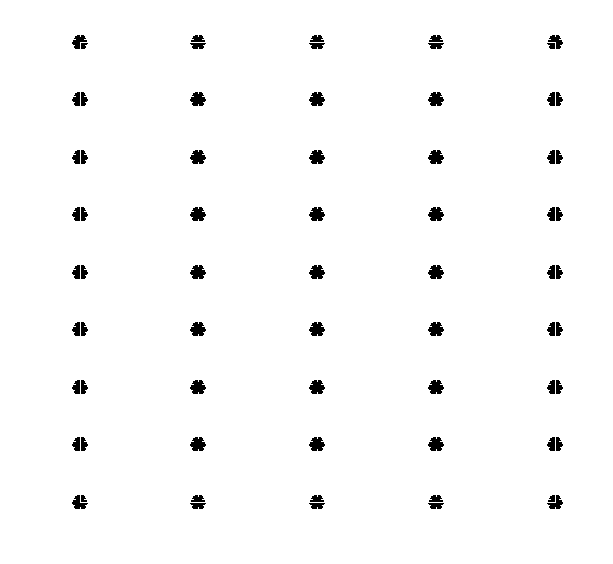} 
     & $\ \cup$ & \includegraphics[width=1.7cm]{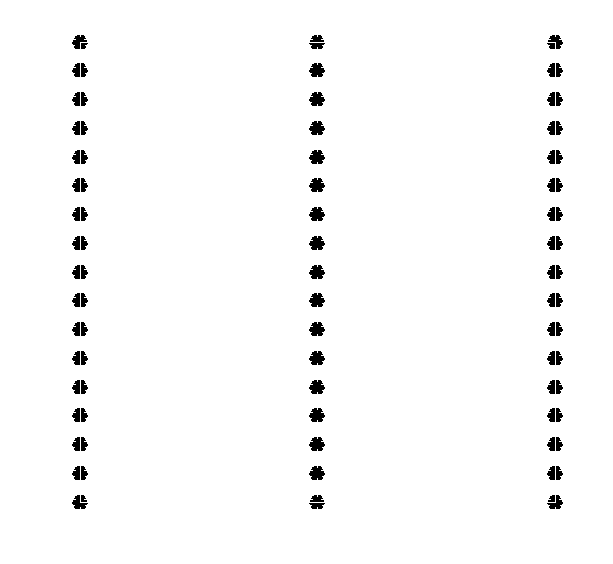}  
\end{tabular}
\end{center}
\caption{Sparse grid $\Upsilon_{4,2}$ via (3). } \label{Fig:Sparse grid def}
\end{figure}

As can easily be seen, the sparse grid treatment itself already requires the use of anisotropic basis functions. For this reason, we will use anisotropic Gaussian functions for Q-SIK. Let $A_{{\bf l}}  = {\rm diag} \, (2^{l_1},2^{l_2},\cdots,2^{l_d})$. Then define
$$
\mu_{{\textbf{l}}}({\bf x}) = (2 \pi)^{-d/2} \exp(-\| A_{{\bf l}} {\bf x} \|^2 ),
$$
and $\mu = \mu_I$, where $I$ is the identity matrix.

In accordance with this purpose, we improve on this via a new Q-SIK interpolation formula on a sparse grid for a function of several variables. In accordance with this purpose, we solve a number of sub-anisotropic quasi-interpolation problems on well constructed subgrids and then obtain the quasi-SIK interpolant by combining the resultant sub-quasi interpolants linearly.

In a similar way, all the subgrids for each level can be constructed as per the previous section. That is, we have $\Upsilon_{\textbf{l}}$ with a mesh size $h_{\textbf{l}}$. Then, in order to compute each subgrid quasi-interpolation problem, we need to construct the anisotropic subgrid quasi-interpolant
\begin{equation}
Q_{\textbf{l}}u(\textbf{x})=\sum_{\textbf{z}\in \Upsilon_{\textbf{l}}}u(\textbf{z})\mu_{\textbf{l}} \left(\textbf{x} -\textbf{z}\right),  \quad \textbf{x}\in \mathbb{R}^n .
\end{equation} 
Then we obtain quasi-SIK (Q-SIK) approximation by the combination technique:
\begin{equation}
Q_{n,d}u(\mathbf{x})=\sum_{q=0}^{d-1}(-1)^q {d-1 \choose q} \sum_{|\mathbf{l}|_1=n+(d-1)-q}Q_{{\textbf{l}}}u(\textbf{x}).
\end{equation}
For example, in 2-D the Q-SIK approximation is
\begin{equation}
Q_{n,2} u(\mathbf{x})=\sum_{|\mathbf{l}|_1=n+1}Q_{{\textbf{l}}}u(\textbf{x})-\sum_{|\mathbf{l}|_1=n}Q_{{\textbf{l}}}u(\textbf{x}).
\end{equation}
\begin{figure}
\begin{center}
\begin{tabular}{m{1.8cm}m{0.2cm}m{1.8cm}m{0.2cm}m{1.8cm}m{0.2cm}m{1.8cm}m{0.2cm}m{1.8cm}}
   \includegraphics[width=1.9cm]{gridS4} &=
      &\includegraphics[width=1.8cm]{gridX41} 
     & $\oplus$  &\includegraphics[width=1.8cm]{gridX32} 
     & $\oplus$ & \includegraphics[width=1.8cm]{gridX23} 
     & $\oplus$ & \includegraphics[width=1.8cm]{gridX14}  \\
     & $\ominus$  &\includegraphics[width=1.8cm]{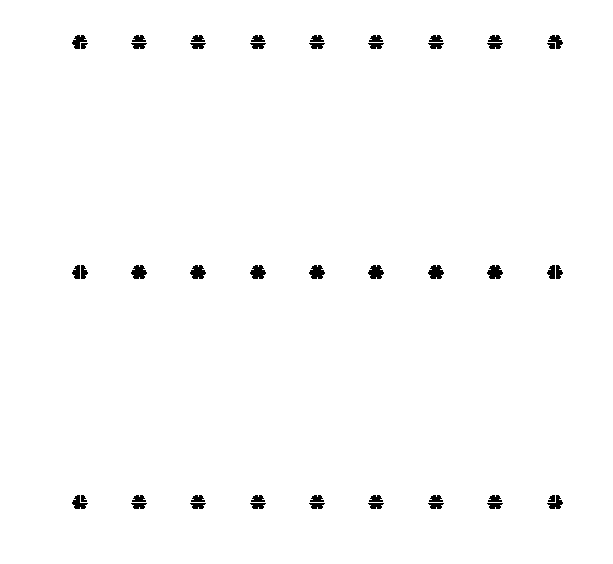} 
     & $\ominus$ & \includegraphics[width=1.8cm]{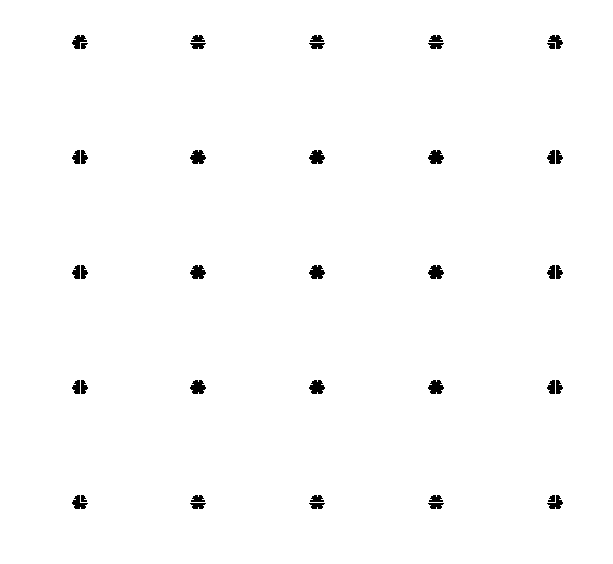} 
     & $\ominus$ & \includegraphics[width=1.8cm]{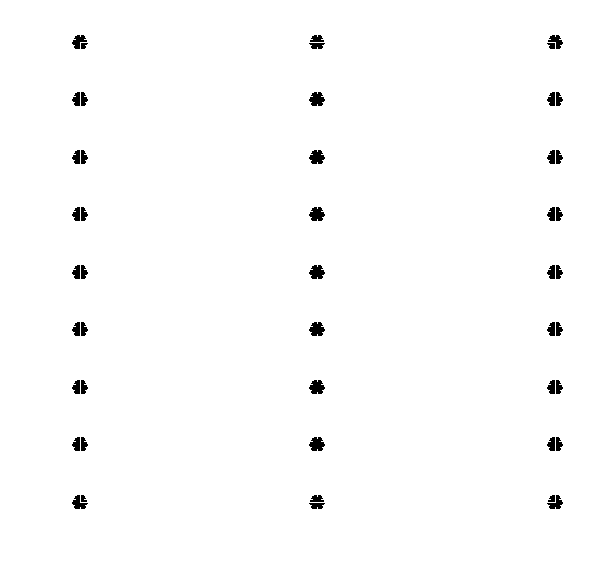}  &
\end{tabular}
\end{center}
\caption{The construction of $\Upsilon_{4,2}$ as an algebraic sum.} \label{Fig:Sparse grid def2}
\end{figure}
The key idea here is that we view the sparse grid in Figure~\ref{Fig:Sparse grid def} as an algebraic sum of grids, as in Figure~\ref{Fig:Sparse grid def2}. We can summarize the algorithm of sparse kernel-based interpolation as  follows.

\begin{algorithm}[H]  
\floatname{algorithm}{Algorithm 1}  
\renewcommand{\thealgorithm}{}  
\caption{Sparse quasi-interpolation with kernels} 
\label{protocol1}  
\begin{algorithmic}[1] 
\REQUIRE Sparse grid data $\lbrace (\textbf{x}_i,u_i),         	u_i=u(\textbf{x}_i), i=1,\cdots N \rbrace$
\STATE Create the subgrids for each level $\textbf{l}\in \mathbb{N}^d$, $\vert \textbf{l}\vert _1 = n,\ldots,n+(d-1)$ and $\Upsilon_{\textbf{l}}$ with mesh size $h_{\textbf{l}}$.
\STATE Compute the anisotropic subgrid quasi-interpolation problems \\ 
$Q_{\textbf{l}}u(\textbf{x})=\sum_{\textbf{z}\in \Upsilon_{\textbf{l}}}u(\textbf{z})\mu_{\textbf{l}} \left(\textbf{x} -\textbf{z}\right).$
\STATE Combine the all subgrid quasi-interpolation problems obtained above\\ $Q_{n,d}u(\mathbf{x})=\sum_{q=0}^{d-1}(-1)^q {d-1 \choose q} \sum_{|\mathbf{l}|_1=n+(d-1)-q}Q_{{\textbf{l}}}u(\textbf{x})$.
\RETURN The Q-SIK approximation $Q_{n,d}u(\mathbf{x})$
\end{algorithmic}  
\end{algorithm} 

As can be seen from the above explanations, the Q-SIK method is very amenable to parallel computation since each sub-interpolation problem can be solved independently. Therefore this method can be applied in a computational cluster or distributed across workstations. Thus, the Q-SIK method provides us with computationally cheap approximation, especially for multidimensional problems.

\section{Multilevel sparse quasi-interpolation with kernels}

As will be seen below (and is observed in SIK), Q-SIK does not converge, even for the function $u \equiv 1$. In the case of the infinite grid in one dimension, this can be seen easily. For a $1/n$ spaced grid the quasi-interpolant is
\beqns
Q_n(x)& = & \sum_{l=-\infty}^\infty \mu (nx-l).
\eeqns
This function is $1/n$-periodic and even and so has a Fourier series
\beqns
Q_n(x) & = & \sum_{k=0}^\infty a_k \cos(2 \pi k n x),
\eeqns
where, for $k \ge 1$,
\beqns
a_k & = & n \int_{0}^{1/n} Q_n(x) \cos(2 \pi k n x) dx \\
& = & n \int_{0}^{1/n} \left ( \sum_{l=-\infty}^\infty \mu (nx-l) \right ) \cos(2 \pi k n x) dx \\
& = & \int_{0}^{1} \left ( \sum_{l=-\infty}^\infty \mu (y-l) \right ) \cos(2 \pi k y) dy \\
& = & \int_{-\infty}^\infty \mu(y) \cos(2 \pi k y) \\
& = & 2\exp(-2 \pi^2 k^2),
\eeqns
the Fourier coefficient of the normal distribution. Additionally $a_0=1$. Hence, for any $n$ 
\beqns
Q_n(0) & = & 1+2\sum_{k=0}^\infty \exp(-2 \pi^2 k^2),
\eeqns
so that there is no convergence as $n \rightarrow \infty$.

Thus we adopt a multilevel refinement strategy, which we call Q-MuSIK. The only difference between MuSIK and Q-MuSIK is that we will use the Q-SIK method instead of the SIK method for each subgrid approximation problem.

\begin{figure}[t!]
  \centering
  \subfloat[$\mathcal{S}_{1,2}$]{\label{figur:1}\includegraphics[width=30mm]{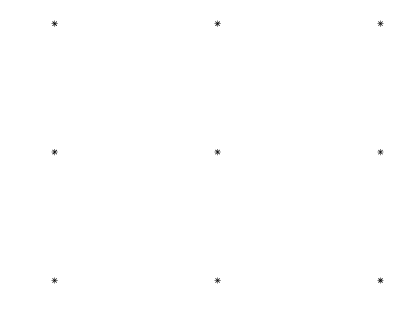}}
  \subfloat[$\mathcal{S}_{2,2}$]{\label{figur:2}\includegraphics[width=30mm]{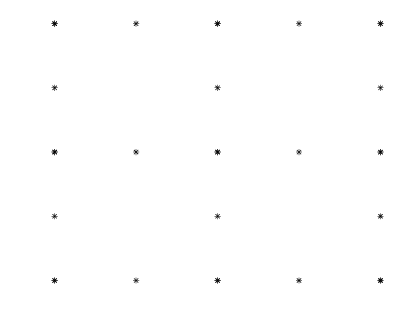}}
  \\
  \subfloat[$\mathcal{S}_{3,2}$]{\label{figur:3}\includegraphics[width=30mm]{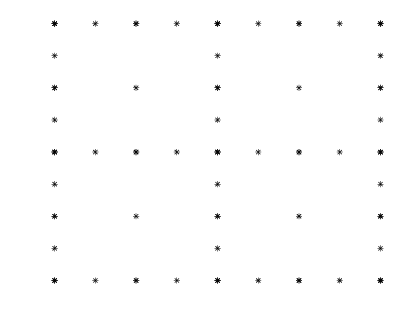}}
  \subfloat[$\mathcal{S}_{4,2}$]{\label{figur:4}\includegraphics[width=30mm]{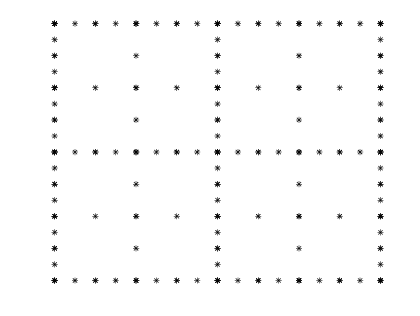}}
  \\
  \subfloat[$\mathcal{S}_{5,2}$]{\label{figur:5}\includegraphics[width=30mm]{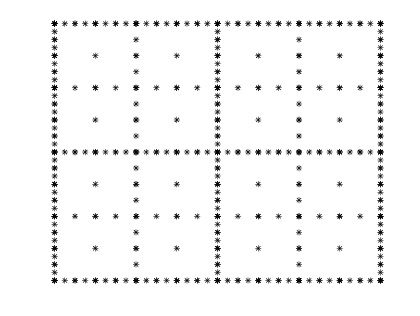}}
  \subfloat[$\mathcal{S}_{6,2}$]{\label{figur:6}\includegraphics[width=30mm]{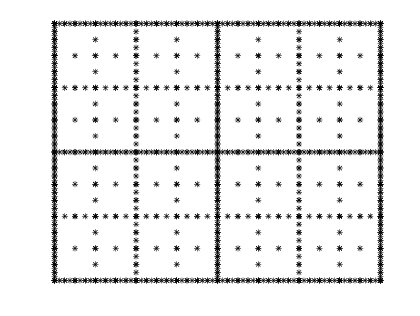}}
  \\
  \subfloat[$\mathcal{S}_{7,2}$]{\label{figur:7}\includegraphics[width=30mm]{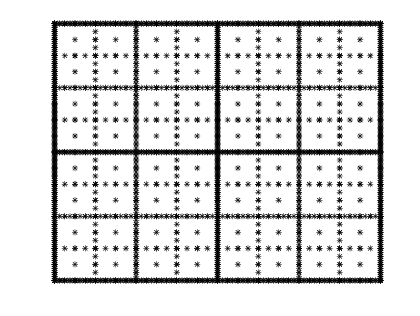}}
  \subfloat[$\mathcal{S}_{8,2}$]{\label{figur:8}\includegraphics[width=30mm]{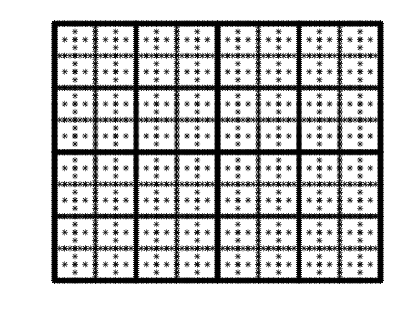}}
\caption{Sparse grid decomposition for levels 1 to 8 in 2 dimensions.} \label{figgrids}
\end{figure}

In contrast to the single level quasi-sparse interpolation with kernels method discussed in the previous chapter, we will now use the nested sparse grids which sort from the lower level to the higher level. In other words, sparse grids with increasingly greater data densities provide us with a  hierarchical decomposition of the sparse grid data,
\begin{equation*}
\Upsilon_{i,d}\subset \Upsilon_{i+1,d}, \quad i=1,2,\cdots.
\end{equation*}
The hierarchical decomposition in two dimensions can be seen in Figure~\ref{figgrids}.

After obtaining the sparse grid decomposition, the sparse grid interpolation $Q_{n_0,d}$ needs to be evaluated for level $1$. In other words, the coarsest level approximation is $S_{0,d}u=Q_{n_0,d}u$. The first level approximation Then $\Delta_j$ needs to be calculated for each level using Q-SIK on the residual $\Delta_j=Q_{n_0+j,d}(u-S_{j-1,d}u)$ on $\Upsilon_{n_0+j,d}$, and $S_{j,d}u=S_{j-1,d} u + \Delta_j$, for $1\leq j\leq n$. This multilevel iterative refinement algorithm is called Q-MuSIK.

\begin{algorithm}[H]  
\floatname{algorithm}{Algorithm 2}  
\renewcommand{\thealgorithm}{}  
\caption{Quasi multilevel sparse interpolation with kernels} 
\label{protocol1}  
\begin{algorithmic}[1] 
\REQUIRE Sparse grid data decomposition and function $u$
\STATE Initialize the first interpolation value at zero, that is $ S_{0,d}u = Q_{n_0,d}u $.
\STATE Construct the nested sparse grids as $ \Upsilon_{n_0,d}\subset \Upsilon_{n_0+1,d}\subset \cdots \subset \Upsilon_{n_0+n,d}  $.  
\STATE For every value of $ j = 1, 2, \cdots ,n   $,
\begin{itemize}
        \item{Solve $ \Delta_j = Q_{n_{0}+j,d}(u-S_{j-1,d}u) $ on $ \Upsilon_{n_0+j,d} $}
        \item{Update $ S_{j,d}u = S_{j-1,d}u + \Delta_j $}
      \end{itemize}     
\RETURN $S_{n,d}$
\end{algorithmic}  
\end{algorithm}   

\section{Q-MuSIK Quadrature}
Now, we introduce a new quadrature formula by integrating the Q-MuSIK approximation over the unit cube:
$$
I_{\textbf{l}} u(\textbf{x}) = \sum_{\textbf{z}\in \Upsilon_{\textbf{l}}}u(\textbf{z})\int_{[0,1]^d} \mu_{\textbf{l}} \left(\textbf{x} -\textbf{z}\right)d\textbf{x}.
$$
In the case when the Gaussians are tensor products,namely, $\mu_{\textbf{l}} \left(\textbf{x} -\textbf{z}\right)=\mu_{l_1}(x_1-z_{l_1})\times \mu_{l_2}(x_2-z_{l_2})\times \cdots \times \mu_{l_d}(x_d-z_{l_d})  $, we can
straightforwardly compute the weights as follows:
$$
\int_{[0,1]^d} \mu_{\textbf{l}} \left(\textbf{x} -\textbf{z}\right)d\textbf{x}=\prod_{i=1}^{d}\int_0^1\mu_{l_i}(x_i-z_{l_i})dx_i
$$
The integral of a univariate Gaussian is of the form
$$
\phi_l(z_l) = \int_0^1 \mu_l(x-z_l)dx=\int_0^1 (2\pi)^{-1/2}\exp(-( A_{l} (x-z_l)) ^2 )dx,
$$
for some $ A_{l} $, $ m_l \in \mathbb{R}$, which, upon the change of variables $ \xi=A_{l} (x-z_l) $, gives
$$
\phi_l(z_l) =\frac{1}{\sqrt{\pi}}\int_{-A_lz_l}^{A_l(1-z_l)}\exp(-\xi^2)d\xi=\frac{1}{2\sqrt{2}}[\erf(A_l(1-z_l))-\erf(-A_lz_l)],
$$
with the error function erf defined by
\begin{equation}
\erf(x)=\frac{2}{\sqrt{\pi}}\int_{0}^{x} e^{-t^2}dt.
\end{equation}
Here the error function can be computed to arbitrary precision.

In a similar way, all the subgrids for each level can be constructed in a similar manner to the method discussed in previous sections. In other words, we will use the grid set $\Upsilon_{\textbf{l}}$ with mesh size $h_{\textbf{l}}$ constructed via sparse grids. In order to then compute each subgrid quadrature problem, we need to use the anisotropic quasi-quadrature technique, which uses sparse grids. In detail, by using the quadrature construction from the previous section, the anisotropic quasi-quadrature formula of $u$ at the grids $\Upsilon_{\textbf{l}}$ can be defined by
$$
I_{\textbf{l}} u(\textbf{x}) = \sum_{\textbf{z}\in \Upsilon_{\textbf{l}}}u(\textbf{z})\prod_{i=1}^d \phi_{l_i}(z_{l_i}).
$$
where \\
$$
\phi_{l_i}(z_{l_i})=\frac{1}{2\sqrt{2}}[\erf(A_{l_i}(1-z_{l_i}))-\erf(-A_{l_i}z_{l_i})]
$$
We then obtain the Q-SIK quadrature formula by combining all sub-quasi-quadrature problems with the combination technique, that is, the Q-SIK quadrature is defined by
$$
I_{n,d}u(\mathbf{x})=\sum_{q=0}^{d-1}(-1)^q {d-1 \choose q} \sum_{|\mathbf{l}|_1=n+(d-1)-q}I_{{\textbf{l}}}u(\textbf{x}).
$$
For example, in 3-D the quasi-SKI is given by
\begin{equation}
I_{n,3}u(\mathbf{x})=\sum_{|\mathbf{l}|_1=n+2}I_lu(\mathbf{x})-2\sum_{|\mathbf{l}|_1=n+1}I_lu(\mathbf{x})+\sum_{|\mathbf{l}|_1=n}I_lu(\mathbf{x}).
\end{equation}

\section{Numerical Experiments}
In this section, we present a collection of numerical experiments for approximation for $d = 2, 3, 4$, where the implementation of Q-SIK and Q-MuSIK algorithms is assessed and compared against both the standard quasi interpolation on uniform full grids and its standard multilevel version (ML Quasi Interpolation). It is worth remarking that in higher dimensions it becomes harder to compute approximation errors. Thus, for higher dimension we move to assessing the errors in quadrature, which is the subject of the final examples, where we consider $d=4,5,10$. We compare Q-MuSIK with MuSIK for these examples.

\subsection{2-$D$ Example}
The Q-SIK and Q-MuSIK algorithm with Gaussian basis functions is
applied to the test function $P^{2d}(\textbf{x}):[0,1]^2\rightarrow\mathbb{R}$, with
\begin{equation}
  \begin{aligned}
  P^{2d}(\textbf{x})&:=\dfrac{1.25+\cos(5.4x_2)}{6+6(3x_1-1)^2}
  \end{aligned}
\end{equation}
\begin{table}[htbp]
    \renewcommand\multirowsetup{\centering}
    \centering
    \begin{tabular}{cccccc}
      \toprule
      \multirow{2}[4]{0.75cm}{\textbf{SGs}} & \multirow{2}[4]{2cm}{\textbf{NoVs}}  & \multicolumn{2}{c}{\textbf{Q-SIK}} & \multicolumn{2}{c}{\textbf{Q-MuSIK}} \\
      \cmidrule{3-6}
          &       &     \textbf{Max Error} & \textbf{Rms Error} & \textbf{Max Error} & \textbf{Rms Error} \\
      \midrule
      \textbf{9} & \textbf{9}  & \textbf{1.48e-1} & \quad \textbf{4.63e-2} &\quad \textbf{1.48e-1} & \quad \textbf{4.63e-2}  \\
      \textbf{21} & \textbf{39}  & \textbf{8.60e-2} &\quad  \textbf{2.08e-2} &\quad \textbf{4.37e-2} &\quad \textbf{1.43e-2}  \\
      \textbf{49} & \textbf{109}  & \textbf{4.92e-2} &\quad  \textbf{9.56e-3} &\quad \textbf{1.61e-2} &\quad \textbf{4.28e-3}  \\
      \textbf{113} & \textbf{271}  & \textbf{3.32e-2} &\quad  \textbf{6.26e-3} &\quad \textbf{7.66e-3} &\quad \textbf{1.31e-3}  \\
      \textbf{257} & \textbf{641}  & \textbf{2.75e-2} & \quad \textbf{5.60e-3} &\quad \textbf{3.26e-3} &\quad \textbf{4.09e-4}  \\
      \textbf{577} & \textbf{1475}  & \textbf{2.63e-2} &\quad  \textbf{5.48e-3} &\quad \textbf{1.33e-3} &\quad \textbf{1.27e-4}  \\
      \textbf{1281} & \textbf{3333}  & \textbf{2.59e-2} &\quad  \textbf{5.45e-3} &\quad \textbf{5.57e-4} &\quad \textbf{3.73e-5}  \\
      \textbf{2817} & \textbf{7431}  & \textbf{2.66e-2} &\quad  \textbf{5.43e-3} &\quad \textbf{1.77e-4} &\quad \textbf{1.01e-5}  \\
      \textbf{6145} & \textbf{16393}  & \textbf{2.61e-2} &\quad  \textbf{5.42e-3} &\quad \textbf{4.77e-5} &\quad \textbf{2.88e-6}  \\
      
      \bottomrule
    \end{tabular}\vspace*{4mm}
     \caption{Q-SIK and Q-MuSIK results using Gaussian basis functions with shape parameter $ d=0.4 $, test function $P^{2d}(\textbf{x})$, on an equally spaced $160\times 160$ evaluation grid in 2 dimension.}
    \label{tab:niveis}
\end{table}
   
\begin{figure}[h!]
\centering
\begin{tabular}{@{}cc@{}}
\vcenteredhbox{\includegraphics[width=0.78\textwidth]{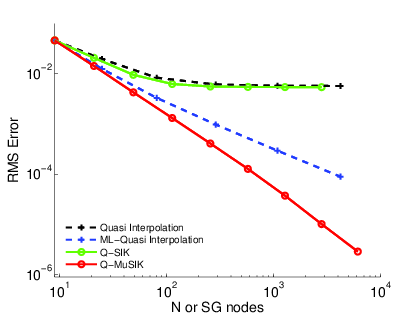}} 
 \end{tabular}
\caption{RMS error versus N(Quasi-interpolation, Multilevel Quasi interpolation) or SG (Q-SIK, Q-MuSIK) nodes using Gaussian basis functions
with $ \rho=0.4 $ Quasi-interpolation (black), Q-SIK (green), Multilevel
Quasi-interpolation (blue) and Q-MuSIK (red) on a $ 160\times 160 $ uniform grid in 2 dimension.}
\label{fig:fig1}
\end{figure}  
\begin{figure}[h!]
\centering
\begin{tabular}{@{}cc@{}}
\vcenteredhbox{\includegraphics[width=0.78\textwidth]{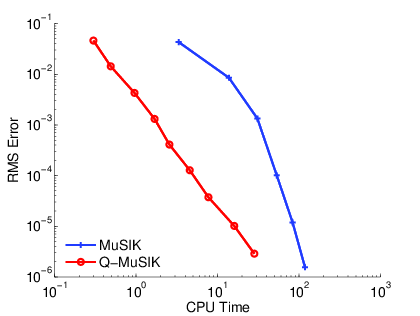}} 
 \end{tabular}
\caption{RMS error versus computational time using Gaussian basis functions with standard deviation $0.4$, MuSIK (blue) and Q-MuSIK (red) on a
$ 160\times 160 $ uniform grid in 2 dimension.}
\label{fig:fig2}
\end{figure}  
 
Here “SGs” stands for the number of sparse grid centers used, “NoVs” represents the total number of nodes visited in the Q-SIK and Q-MuSIK, “Maxerror” and “RMS-error’ are the maximum norm and root-mean-square ($L^2$-norm) errors, respectively, evaluated on a $160 \times 160$ uniform grid. According to the results of Table \ref{tab:niveis} we can see that, as expected, Q-SIK does not converge. In other words, the error of Q-SIK method remain stable after a few interpolation levels. However we see in the same table that Q-MuSIK converges.

In Figure \ref{fig:fig1} we present the comparison four interpolation results which are quasi interpolation, multilevel quasi interpolation, Q-SIK and Q-MuSIK. Obviously, quasi interpolation and Q-SIK do not converge to the target functions. However their multilevel versions shows good convergence behaviour. 

In Figure \ref{fig:fig2} we compare the MuSIK and Q-MuSIK results with regard to computational time. According to this figure Q-MuSIK reaches the same error level as MuSIK in a shorter time. This is because quasi interpolation does not require the solution of any large algebraic system as does an interpolatory scheme such as MuSIK.

\subsection{3-$D$ Example}
\begin{figure}[h!]
\centering
\begin{tabular}{@{}cc@{}}
\vcenteredhbox{\includegraphics[width=0.78\textwidth]{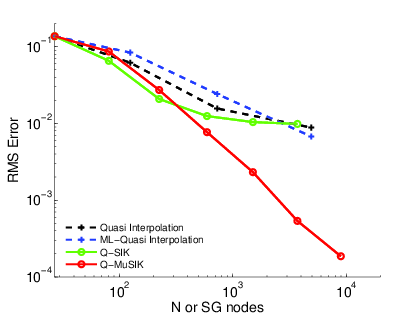}} 
 \end{tabular}
\caption{RMS error versus N (Quasi-interpolation, Multilevel Quasi interpolation) or SG (Q-SIK, Q-MuSIK) nodes using Gaussian basis functions
with standard deviation $ 0.4 $, Quasi-interpolation (black), Q-SIK (green), Multilevel
Quasi-interpolation (blue) and Q-MuSIK (red) on a $ 50\times 50\times 50 $ uniform grid in 3 dimension.}
\label{fig:fig3}
\end{figure}

\begin{figure}[h!]
\centering
\begin{tabular}{@{}cc@{}}
\vcenteredhbox{\includegraphics[width=0.78\textwidth]{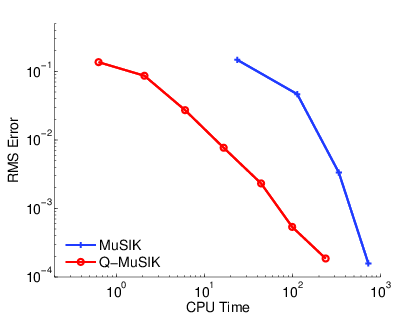}} 
 \end{tabular}
\caption{RMS error versus computational time using Gaussian basis functions
with standard deviation 0.4, MuSIK (blue) and Q-MuSIK (red) on a
$ 50\times 50\times 50 $ uniform grid in 3 dimension.}
\label{fig:fig4}
\end{figure}

The following example is in 3 dimensions evaluated at an equally spaced $50\times 50 \times 50$ evaluation grid.  We consider the interpolation problem of the another test function $G^{3d}(\textbf{x}):[0,1]^3\rightarrow\mathbb{R}$, with 
\begin{equation}
G^{3d}(\textbf{x}):=\dfrac{18}{\pi}e^{-(x_1^2+81x_2^2+x_3^2)}
\end{equation}
Experimental results in 3 dimension indicate that Q-MuSIK shows a good approximation behaviour as is the case with previous example. The corresponding results to those in 2 dimensions can be seen in Figure \ref{fig:fig3} and \ref{fig:fig4}.   

\subsection{4-$D$ Example}
\begin{figure}[h!]
\centering
\begin{tabular}{@{}cc@{}}
\vcenteredhbox{\includegraphics[width=0.78\textwidth]{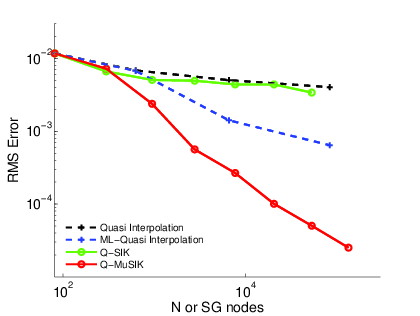}} 
 \end{tabular}
\caption{RMS error versus N(Quasi-interpolation, Multilevel Quasi-interpolation) or SG (Q-SIK, Q-MuSIK) using Gaussian basis functions with standard deviation 0.4: Quasi-interpolation (black), Q-SIK (green), Multilevel Quasi-interpolation (blue) and Q-MuSIK (red) on a $21\times 21\times 21 \times 21$ uniform grid.}
\label{fig:fig5}
\end{figure}

Let us now consider 4 dimension non-tensor product test functions $ H^{4d}(\textbf{x}): [0,1]^4\rightarrow \mathbb{R} $, with
\begin{equation}
H^{4d}(\textbf{x})=\sin(x_1^2x_2^2x_3^2x_4^2).
\end{equation} 
In Figure \ref{fig:fig5} the root mean-square error of Quasi Interpolation, ML Quasi Interpolation, Q-SIK and Q-MuSIK, respectively, for $ H^{4d}(\textbf{x})$, are plotted against the number of data sites. We see similar performance as in 2 and 3 dimensions. 

\subsection{Quadrature Example}
In this section we give a number of examples for quadrature. We use two tensor product examples, respectively in 5 and 10 dimensions, and two non-tensor product examples. We consider the following non-tensor product example
\begin{equation}
M^{3d(\textbf{x})}= \sin(x_1x_2x_3).
\end{equation}

\begin{figure}[h!]
\centering
\begin{tabular}{@{}cc@{}}
\vcenteredhbox{\includegraphics[width=0.78\textwidth]{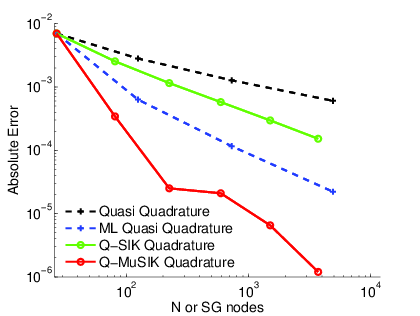}} 
 \end{tabular}
\caption{Absolute error versus N(Quasi-interpolation, Multilevel Quasi interpolation) or SG (Q-SIK, Q-MuSIK) nodes using Gaussian basis functions with standard deviation 0.4 on the domain $ [0,1]^3 $.}
\label{fig:fig6}
\end{figure}

In Figure \ref{fig:fig6} we give the Quasi, Multilevel Quasi, Q-SIK and Q-MuSIK quadrature results for $ M^{3d(\textbf{x})} $ with regard to the corresponding exact values of the integral of $ M^{3d(\textbf{x})} $ on the domain $[0,1]^{3}$, which is (with 16 digits accuracy) $  0.122434028745371 $. 

\begin{figure}[htp]
  \centering
  \subfloat[$F^{4d}$]{\label{figur:1}\includegraphics[width=70mm]{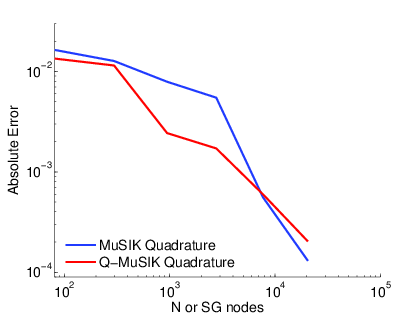}}
  \subfloat[$F^{4d}$]{\label{figur:2}\includegraphics[width=70mm]{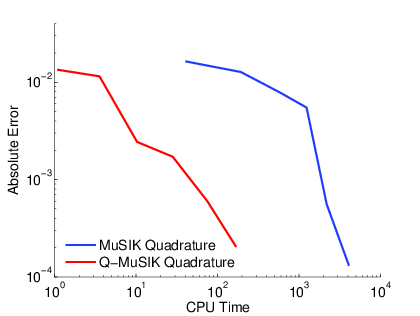}}
  \\
  \subfloat[$T^{5d}$]{\label{figur:3}\includegraphics[width=70mm]{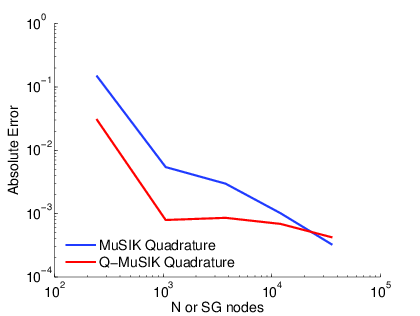}}
  \subfloat[$T^{5d}$]{\label{figur:4}\includegraphics[width=70mm]{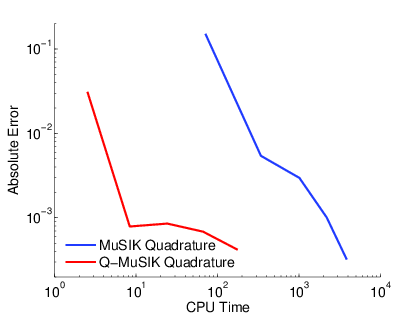}}
  \\
  \subfloat[$K^{10d}$]{\label{figur:5}\includegraphics[width=70mm]{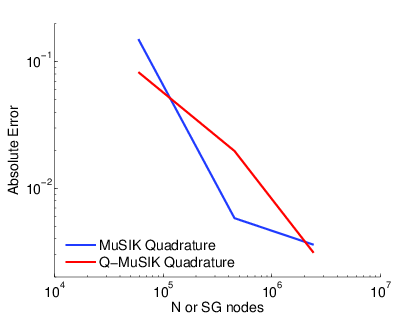}}
  \subfloat[$K^{10d}$]{\label{figur:6}\includegraphics[width=70mm]{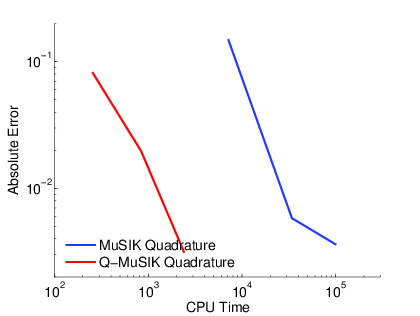}}
\caption{Absolute error versus SG and CPU Time(Q-MuSIK quadrature, Q-MuSIK quadrature) nodes using Gaussian basis functions with $\rho=0.4$ and $c=0.45$: MuSIK quadrature (blue) and Q-MuSIK quadrature (red) on the domain $[0,1]^{4,5,10}$.}
\label{figur10}
\end{figure}

In Figure \ref{figur10} we compare the MuSIK quadrature and Q-MuSIK quadrature in terms of absolute error versus degree of freedom and absolute error versus evaluation time. As discussed in the previous examples Q-MuSIK has better performance especially at the low level of computations. In these comparison we have used four, five and ten dimensional test functions which are given below:
\begin{eqnarray*}
F^{4d(\textbf{x})}&=&\frac{3}{4} e^{-(9x_1-2)^2-(9x_2-2)^2-(9x_3-2)^2/4-(9x_1-2)^2/8}\\
&+&\frac{3}{4} e^{-(9x_1+1)^2/49-(9x_2+1)^2/10-(9x_3+1)^2/29-(9x_1+1)^2/39}\\
&+&\frac{1}{2} e^{-(9x_1-7)^2/4-(9x_2-3)^2-(9x_3-5)^2/2-(9x_1-5)^2/4} \\
&-&\frac{1}{5} e^{-(9x_1-4)^2/4-(9x_2-7)^2-(9x_3-5)^2-(9x_1-5)^2}\\\\,
T^{5d}(\textbf{x})&=& \prod_{i=1}^{4} \max(x_i-1/2)\\\\,
K^{10d}(\textbf{x})&=&\prod_{i=1}^{10}e^{-x_i(1-x_i)},
\end{eqnarray*}  
and the corresponding exact integral values of these functions on the domain $[0,1]^{4}$, $[0,1]^{5}$ and $[0,1]^{10}$ with 8 digits accuracy are $0.07766696$, $0.001953125$ and $0.19427907$ respectively. These results confirms that the Q-MuSIK method provide us more rapid results for the same amount of error as the comparitor methods.

\section{Concluding remarks}

A new quasi multilevel sparse interpolation with Gaussian kernel for, particularly high-dimensional, both interpolation and numerical integration step is recommended and examined. One of the most significant positive aspects of quasi-interpolation is that it can yield a solution directly with no need to solve large algebraic systems. In addition, it has a number of desirable features, such as good approximation behaviour, easy evaluation and computation. The Q-MuSIK algorithm couples these good properties and the complexity of the sparse grid technique by using  a multilevel algorithm. Q-MuSIK appears to be usually superior over classical the quasi interpolation technique in terms of run time and complexity especially in high dimension. In addition to this the Q-MuSIK algorithm computes approximations more quickly than the MuSIK method in problems with the same number of degrees of freedom because of the properties of quasi-interpolation listed above.
\bibliographystyle{siam}
\bibliography{Thesisbibdatanew}

\end{document}